\documentclass[twoside,leqno,a4paper]{article}
\usepackage{amsmath}
\usepackage{amsthm}
\usepackage{amssymb}
\usepackage{amscd}
\usepackage{graphicx}
\usepackage{epsfig}

\usepackage{latexsym}
\usepackage{amsfonts}
\input xy
\xyoption{all} \tolerance=500
=1in \textwidth 15cm \topmargin 5mm\textheight24cm \advance\topmargin by -\headheight\advance\topmargin by
-\headsep

\input xy
\xyoption{all} \tolerance=500

\font\black=cmbx10 \font\sblack=cmbx7 \font\ssblack=cmbx5 \font\blackital=cmmib10  \skewchar\blackital='177
\font\sblackital=cmmib7 \skewchar\sblackital='177 \font\ssblackital=cmmib5 \skewchar\ssblackital='177
\font\sanss=cmss10 \font\ssanss=cmss8 %scaled 900
\font\sssanss=cmss8 scaled 600 \font\blackboard=msbm10 \font\sblackboard=msbm7 \font\ssblackboard=msbm5
\font\caligr=eusm10 \font\scaligr=eusm7 \font\sscaligr=eusm5  \font\fraktur=eufm10
\font\sfraktur=eufm7 \font\ssfraktur=eufm5

\font\bsymb=cmsy10 scaled\magstep2
\def\all#1{\setbox0=\hbox{\lower1.5pt\hbox{\bsymb
       \char"38}}\setbox1=\hbox{$_{#1}$} \box0\lower2pt\box1\;}
\def\exi#1{\setbox0=\hbox{\lower1.5pt\hbox{\bsymb \char"39}}
       \setbox1=\hbox{$_{#1}$} \box0\lower2pt\box1\;}

\def\tx#1{{\fam0\relax#1}}

\newfam\bifam
\textfont\bifam=\blackital \scriptfont\bifam=\sblackital \scriptscriptfont\bifam=\ssblackital
\def\bi#1{{\fam\bifam\relax#1}}

\newfam\blfam
\textfont\blfam=\black \scriptfont\blfam=\sblack \scriptscriptfont\blfam=\ssblack

\newfam\bbfam
\textfont\bbfam=\blackboard \scriptfont\bbfam=\sblackboard \scriptscriptfont\bbfam=\ssblackboard

\newfam\ssfam
\textfont\ssfam=\sanss \scriptfont\ssfam=\ssanss \scriptscriptfont\ssfam=\sssanss
\def\sss#1{{\fam\ssfam\relax#1}}

\newfam\clfam
\textfont\clfam=\caligr \scriptfont\clfam=\scaligr \scriptscriptfont\clfam=\sscaligr

\newfam\frfam
\textfont\frfam=\fraktur \scriptfont\frfam=\sfraktur \scriptscriptfont\frfam=\ssfraktur

\def\hpb#1{\setbox0=\hbox{${#1}$}
    \copy0 \kern-\wd0 \kern.2pt \box0}
\def\vpb#1{\setbox0=\hbox{${#1}$}
    \copy0 \kern-\wd0 \raise.08pt \box0}

\def\pmb#1{\setbox0\hbox{${#1}$} \copy0 \kern-\wd0 \kern.2pt \box0}
\def\pmbb#1{\setbox0\hbox{${#1}$} \copy0 \kern-\wd0
      \kern.2pt \copy0 \kern-\wd0 \kern.2pt \box0}
\def\pmbbb#1{\setbox0\hbox{${#1}$} \copy0 \kern-\wd0
      \kern.2pt \copy0 \kern-\wd0 \kern.2pt
    \copy0 \kern-\wd0 \kern.2pt \box0}
\def\pmxb#1{\setbox0\hbox{${#1}$} \copy0 \kern-\wd0
      \kern.2pt \copy0 \kern-\wd0 \kern.2pt
      \copy0 \kern-\wd0 \kern.2pt \copy0 \kern-\wd0 \kern.2pt \box0}
\def\pmxbb#1{\setbox0\hbox{${#1}$} \copy0 \kern-\wd0 \kern.2pt
      \copy0 \kern-\wd0 \kern.2pt
      \copy0 \kern-\wd0 \kern.2pt \copy0 \kern-\wd0 \kern.2pt
      \copy0 \kern-\wd0 \kern.2pt \box0}

\mathchardef\za="710B  %\alpha
\mathchardef\zb="710C  %\beta
\mathchardef\zg="710D  %\gamma
\mathchardef\zd="710E  %\delta
\mathchardef\zve="710F %\epsilon
\mathchardef\zz="7110  %\zeta
\mathchardef\zh="7111  %\eta
\mathchardef\zvy="7112 %\theta
\mathchardef\zi="7113  %\iota
\mathchardef\zk="7114  %\kappa
\mathchardef\zl="7115  %\lambda
\mathchardef\zm="7116  %\mu
\mathchardef\zn="7117  %\nu
\mathchardef\zx="7118  %\xi
\mathchardef\zp="7119  %\pi
\mathchardef\zr="711A  %\rho
\mathchardef\zs="711B  %\sigma
\mathchardef\zt="711C  %\tau
\mathchardef\zu="711D  %\upsilon
\mathchardef\zvf="711E %\phi
\mathchardef\zq="711F  %\chi
\mathchardef\zc="7120  %\psi
\mathchardef\zw="7121  %\omega
\mathchardef\ze="7122  %\varepsilon
\mathchardef\zy="7123  %\vartheta
\mathchardef\zf="7124  %\varomega
\mathchardef\zvr="7125 %\varrho
\mathchardef\zvs="7126 %\varsigma
\mathchardef\zf="7127  %\varphi
\mathchardef\zG="7000  %\Gamma
\mathchardef\zD="7001  %\Delta
\mathchardef\zY="7002  %\Theta
\mathchardef\zL="7003  %\Lambda
\mathchardef\zX="7004  %\Xi
\mathchardef\zP="7005  %\Pi
\mathchardef\zS="7006  %\Sigma
\mathchardef\zU="7007  %\Upsilon
\mathchardef\zF="7008  %\Phi
\mathchardef\zW="700A  %\Omega
\mathchardef\zC="7009  %\Psi

\newcommand{\be}{\begin{equation}}
\newcommand{\ee}{\end{equation}}

\newcommand{\bea}{\begin{eqnarray}}
\newcommand{\eea}{\end{eqnarray}}
\newcommand{\beas}{\begin{eqnarray*}}
\newcommand{\eeas}{\end{eqnarray*}}
\def\*{{\textstyle *}}
\newcommand{\R}{{\mathbb R}}

\newcommand{\GL}{\mathrm{GL}}

\newcommand{\we}{\wedge}

\newcommand{\s}{{\textstyle *}}

\newcommand{\pa}{\partial}
\newcommand{\ti}{\times}

\def\bd{{\bi d}}

\def\sT{{\sss T}}

\def\xi{\tx{i}}

\newdir{ (}{{}*!/-5pt/@^{(}}

\def\s*{{\scriptstyle *}}

\newcommand{\bfr}{\begin{frame}}
\newcommand{\efr}{\end{frame} }

\def\wu{{\operatorname{deg}}}

\newdir{|>}{%
!/4.5pt/@{|}*:(1,-.2)@^{>}*:(1,+.2)@_{>}}
\newdir{ (}{{}*!/-5pt/@^{(}}

\def\Linr{{\operatorname{L}}} %full linearisation functor
\def\pLinr{{\operatorname{l}}} %partial linearisation functor

\newcommand{\SymmVB}{\operatorname{SymmVB}}
\newcommand{\GrB}{\operatorname{GrB}}

\newcommand{\VB}{\operatorname{VB}}
\newcommand{\GrL}{\operatorname{GrL}}

\newcommand{\qDA}{{q^D_A}}
\newcommand{\qDB}{{q^D_B}}

\newcommand{\qA}{{q_A}}
\newcommand{\qB}{{q_B}}

\newcommand{\Aplus}{~\ {+}_{A}\ ~}

\newcommand{\Atimes}{\ {\cdot}_{A}\ }

\newcommand{\zeroDA}{\tilde{0}^A}

\newcommand{\Bplus}{~\ {+}_{B}\ ~}

\newcommand{\Btimes}{\ {\cdot}_{B}\ }
\newcommand{\zeroDB}{\tilde{0}^B}
\newcommand{\lon }{\,\rightarrow\,}

\numberwithin{equation}{section}

\newtheorem{Theorem}{Theorem}[section]
\newtheorem{Corollary}{Corollary}
\theoremstyle{definition}
\newtheorem{Definition}[Theorem]{Definition}

\newtheorem{Example}[Theorem]{Example}

\begin{document}
\bibliographystyle{plain}

\author{\\
        Andrew James Bruce$^1$\\ Katarzyna  Grabowska$^2$\\
        Janusz Grabowski$^1$\\
        \\
         $^1$ {\it Institute of Mathematics,
                Polish Academy of Sciences }\\
         $^2$ {\it Faculty of Physics, University of Warsaw} }

\date{\today}
\title{Introduction to graded bundles\thanks{Research funded by the  Polish National Science Centre grant under the contract number DEC-2012/06/A/ST1/00256.  }}
\maketitle
\begin{abstract}
We review the concept of a  graded bundle as a natural generalisation of a vector bundle. Such geometries are   particularly nice examples of more general graded manifolds. With hindsight there are many examples of graded bundles that appear in the existing literature.  We  start with a discussion of graded spaces, passing through graded bundles and their linearisation to end with weighted (graded) Lie groupoids/algebroids.
\end{abstract}

\bigskip
\begin{small}

\smallskip

\end{small}

\tableofcontents

\section{Introduction}
In recent years there has been an interest in various notions of `graded geometry', largely this interest has been inspired by the BV-BRST and BFV formalisms of gauge field theory, as well as the related AKSZ method. From a pure mathematical perspective, graded geometries have been useful in understanding Lie algebroids, Courant algebroids, Lie bialgebroids and related structures.

Our general understanding of graded supergeometry is in the sense of Voronov \cite{Voronov:2001qf}. Loosely, a graded supermanifold is a supermanifold for which the structure sheaf carries in addition to the $\mathbb{Z}_{2}$-grading of Grassmann parity an independent weight from $\mathbb{Z}$. That is, we can always find homogeneous local coordinates that are assigned both a Grassmann parity and a weight, and changes of coordinates respect these assignments.

Among the  possible graded geometries, (super)manifolds that have a non-negative grading seem to play a prominent r\^{o}le in both the general theory and applications of graded manifolds. In this review we will stick to a  class of non-negatively graded manifolds that have particularly nice properties:  we refer to these manifolds as \emph{graded bundles}  (cf. \cite{Grabowski:2012}).

The key concept in the theory of graded bundles is that of a \emph{homogeneity structure}. Under the assumption that our non-zero weight coordinates can take on any real value, we see that a graded structure leads to an action of the multiplicative monoid of reals: this is just a `dilation' defined by the weight of the coordinates. For example, the higher order tangent bundles of a manifold $\sT^{k} M$, i.e. the higher jets of curves in $M$, come with a canonical graded structure of this kind. It is vital that we consider non-negative weights here, otherwise  we cannot think in terms of the action of the monoid of multiplicative reals.

Graded bundles, which we will define carefully in due course, are graded manifolds for which the graded structure is inherited from the smooth action of the multiplicative monoid of reals, which is known as a \emph{homogeneity structure} (cf. \cite{Grabowski:2012}). If the homogeneity structure is regular, then following \cite{GR0} we know we in fact have the structure of the total space of a vector bundle.  In this precise sense, graded bundles represent a natural generalisation of the notion of a vector bundle. The example to always keep in mind here is the passing from a tangent bundle to a higher order tangent bundle. There are many other natural examples of graded bundles and some of these will be presented in this review.

Phrasing the theory in terms of a homogeneity structure is not only conceptually neat, but it also allows for a clear understanding of $n$-fold graded bundles in terms of commuting actions. For example, the notion of a double vector bundle is clear, we have a pair of regular homogeneity structures and their actions commute. In other words, the homogeneity structures are compatible. As we shall see, using homogeneity structures leads to very clear and concise notions of weighted Lie groupoids and algebroids as the classical structures equipped with a compatible homogeneity structure. The only thing to decide is a reasonable notion of compatibility. We will see that the compatibility is simply a natural condition that the action of the homogeneity structure be a morphism in the appropriate category, for instance the category  Lie groupoids or Lie algebroids.

We will stick to pure mathematical aspects of the theory of graded bundles throughout this review and we only remark that many of the concepts we present here have already found applications in higher order mechanics, see \cite{Bruce:2014}. It is certainly expected that further applications of graded bundles and related structures will be found. There are many questions that are still open and it is not completely clear how far we can push the philosophy that graded bundles are `higher' or `nonlinear' vector bundles.

This review is  based on the work of the third author in collaboration with the first two authors and Rotkiewicz: \cite{Bruce:2014, Bruce:2014a, Bruce:2016a, Grabowski:2012}. Proofs of the various statements made in this review will not be given, instead we direct the interested reader to the original literature. \\

\noindent \textbf{Arrangement:} In Section \ref{sec:Graded spaces} we introduce the notion of a graded space as a natural generalisation of a vector space.
Graded spaces are then used in Section \ref{sec:GradedBundles} to define the notion of a graded bundle, which we view as a generalisation of a vector bundle.

Double structures, which include the well-known double vector bundles, are the subject of Section \ref{sec:DoubleStructures}.
These double structures appear in Section \ref{sec:Linearisation} where we discuss the linearisation of a graded bundle.

The final part of this review, Section \ref{sec:WeightedLie}, is devoted to the notion of weighted Lie groupoids and algebroids, which are examples of geometrical objects equipped with a compatible homogeneity structure.

\section{Graded spaces}\label{sec:Graded spaces}
According to textbook definition, a real vector space is a set $E$ with a distinguished element $0^E$, equipped with two operations: an addition
$$+:E\ti E\to E\,,\quad (u,v)\mapsto u+v\,,$$
and a multiplication by scalars
$$h:\R\ti E\to E\,,\quad h(t,v)=h_t(v)=t\cdot v=tv\,,$$
satisfying a list of axioms.  For instance, $(E,+)$ is a commutative group with $0^E$ being the neutral element. The \emph{homotheties} $h_t$ satisfy
\begin{equation}\label{gr:1} h_t\circ h_s=h_{ts}\,,\quad\text{and}\quad \forall v\in E\;\; h_0(v)=0^E.\end{equation}

Every finite-dimensional vector space is also a differential manifold, since any basis provides it with global coordinates. Let us notice that to distinguish real vector spaces among differentiable manifolds, a single operation of the above two is enough.

If we know that the addition comes from a real vector space structure, we get the multiplication by natural numbers in the obvious way:
$nv=v+\cdots + v$, and we easily extend it to integers by $(-n)v=n(-v)$. The multiplication by rational numbers, $(m/n)v$ we obtain as the unique solution of the equation $nx=mv$. Assuming differentiability (in fact, continuity) of $h$, we extend this multiplication to all reals uniquely.

If, in turn, we know the map $h$, i.e. multiplication by reals, we use a consequence of Euler's Homogeneous Function Theorem: any differentiable $f:\R^n\to\R$ is homogeneous of degree 1, i.e. $f(t\cdot x)=t\cdot f(x)$, if and only if $f$ is linear.
Thus, from the multiplication by reals on $E$ we get the dual space $E^*$, where the addition is well defined, and consequently the addition on $E=(E^*)^*$.

But how we can recognize that the map $h$ satisfying (\ref{gr:1}) comes from a real vector space structure? A simple answer to this question is given by the following.
\begin{Theorem}[Grabowski-Rotkiewicz \cite{GR0}]\label{theorem:0}
A smooth action $h:\R\ti E\to E$ of the monoid of multiplicative reals on a manifold $E$,  satisfying $h_0(v)=0^E$ for some $0^E\in E$ and all $v\in E$, comes from a real vector space structure on $E$ if and only if $h$ is \emph{regular}, i.e.
$$\frac{\partial h}{\partial t}(0,v)=0\quad \Leftrightarrow\quad v=0^E\,.$$
In this case, the real vector space structure is unique.
\end{Theorem}
Vector spaces can be generalized by accepting wider class of maps $h$ than just the regular ones.

Consider now just a smooth action $h:\R\ti F\to F$ of the monoid $(\R,\cdot)$ on a manifold $F$ and assume that $h_0(F)=0^F$ for some element $0^F\in F$. Such an action we will call a \textit{ homogeneity structure}. The set $F$ with a homogeneity structure will be called a \textit{graded space}. The reason for the name is the following theorem

\begin{Theorem}[Grabowski-Rotkiewicz \cite{Grabowski:2012}]\label{theorem:1}
Any graded space $(F,h)$ is diffeomorphically equivalent to a \emph{dilation structure}, i.e. to a certain $(\R^\bd,h^\bd)$, where
$\bd=(d_1,\dots,d_k)$, with positive integers $d_i$,   and $\R^\bd=\R^{d_1}\ti\cdots\ti\R^{d_k}$
is equipped with the dilation action $h^\bd$ of multiplicative reals given by
$$h^\bd_t(y_1,\dots,y_k)=(t\cdot y_1,\dots, t^k\cdot y_k)\,,\quad y_i\in\R^{d_i}\,.$$
In other words, $F$ can be equipped with a system of (global) coordinates $(y_i^j)$, $i=1\dots, k$, $j=1,\dots, d_i$, such that $y_i^j$ is  {homogeneous of degree $i$} with respect to the homogeneity structure $h$:
$$y_i^j\circ h_t=t^i\cdot y_i^j\,.$$
Of course, in these coordinates $0^F=(0,\dots, 0)$.
\end{Theorem}

It is natural to call a \textit{morphism between graded spaces} $(F_a,h^a)$, $a=1,2$, a smooth map $\Phi:F_1\to F_2$ which intertwines the homogeneity structures:
\begin{equation}\label{gr:2}\Phi\circ h_t^1=h_t^2\circ\Phi
\end{equation}
The $(\R,\cdot)$-action, restricted to positive reals, gives a one-parameter group of diffeomorphisms of $F$, thus is generated by a vector field $\nabla_F$. It is called the \textit{ weight vector field} as it completely determines the weights (degrees) of homogeneous functions.
In homogeneous local coordinates the weight vector field reads
\begin{equation}\label{gr:3}\nabla_F=\sum_{w,j}w\, y_w^j\partial_{y^j_w}\,.\end{equation}
We say that a function $f$ on $F$ is  \textit{homogeneous of degree $w$} (\textit{has weight $w$}) if $\nabla_F(f)=w\cdot f$. A smooth map $\Phi:F_1\to F_2$ is a morphism of graded spaces if and only if it relates the corresponding weight vector fields.

Note that automorphisms of $(\R^\bd,h^\bd)$ need not to be linear, so the category of graded spaces is different from that of vector spaces.   For instance, if $(y,z)\in\R^2$ are coordinates of degrees $1,2$, respectively, then the map
$$\R^2\ni(y,z)\mapsto (y,z+y^2)\in\R^2$$
is an automorphism of the homogeneity structure, but is nonlinear.

\section{Graded bundles}\label{sec:GradedBundles}

A \textit{vector bundle} is a locally trivial fibration $\zt:E\to M$ which, locally over some open subsets $U\subset M$, reads $\zt^{-1}(U)\simeq U\ti\R^n$ and admits an atlas in which local trivialisations transform linearly in fibers:
\begin{equation}\label{bu:1} U\cap V\ti\R^n\ni(x,y)\longmapsto(x,A(x)y)\in U\cap V\ti\R^n\,,\quad
A(x)\in\GL(n,\R).\end{equation}
The latter property can also be expressed in the terms of the gradation in which base coordinates $x$ have degrees $0$, and `linear coordinates' $y$ have degree $1$. Linearity in $y's$ is now equivalent to the fact that changes of coordinates respect the degrees. Morphisms in the category of vector bundles are represented by commutative diagram of smooth maps
\begin{equation}\label{bu:2} \xymatrix{
E_1\ar[rr]^{\Phi} \ar[d]^{\zt_1} && E_2\ar[d]^{{\zt_2}} \\
M_1\ar[rr]^{{\varphi}} && M_2 }
\end{equation}
being linear (homogeneous) in fibres.

A straightforward generalization is the concept of a \textit{graded bundle} $\zt:F\to M$
of rank $\bd$, with a local trivialization by $U\ti\R^\bd$, and with the difference that the transition functions of local trivialisations:
$$U\cap V\ti\R^\bd\ni(x,y)\longmapsto(x,A(x,y))\in U\cap V\ti\R^\bd\,,$$
respect the weights of coordinates $(y^1,\dots,y^{|\bd|})$ in the fibres, i.e. $A(x,\cdot)$ are automorphisms of the graded space $(\R^\bd,h^\bd)$. In other words, a graded bundle of rank $\bd$ is a locally trivial fibration with fibers modelled on the graded space $\R^\bd$.

\begin{Theorem}[Grabowski-Rotkiewicz \cite{Grabowski:2012}] $A(x,y)$ must be polynomial in homogeneous fiber coordinates $y$'s, i.e.
any graded bundle is a \textit{polynomial bundle}.
\end{Theorem}
As these polynomials need not to be linear, graded bundles do not have, in general,
vector space structure in fibers. If all $w_i\le r$, we say that the graded bundle is \emph{of degree} $r$. In the above terminology, {vector bundles are just graded bundles of degree $1$}.

\begin{Example} Consider the second-order tangent bundle $\sT^2M$, i.e. the bundle of second jets of smooth maps $(\R,0)\to M$. Writing Taylor expansions of curves in local coordinates $(x^A)$ on $M$:
\begin{equation*}x^A(t)=x^A(0)+\dot x^A(0)t+\ddot x^A(0)\frac{t^2}{2}+o(t^2)\,,\end{equation*}
we get local coordinates $(x^A,\dot x^B,\ddot x^C)$ on $\sT^2M$, which transform
\beas
x'^A&=&x'^A(x)\,,\\
\dot x'^A&=&\frac{\pa x'^A}{\pa x^B}(x)\,\dot x^B\,,\\
\ddot x'^A&=&\frac{\pa x'^A}{\pa x^B}(x)\,\ddot x^B+\frac{\pa^2 x'^A}{\pa x^B\pa x^C}(x)\,\dot x^B\dot x^C\,.
\eeas
This shows that associating with $(x^A,\dot x^B,\ddot x^C)$ the weights $0,1,2$, respectively, will give us a graded bundle structure of degree $2$ on $\sT^2M$. Note that, due to the quadratic terms above, this is not a vector bundle.
All this can be generalised to higher tangent bundles $\sT^kM$.
\end{Example}

\begin{Example}[\cite{GGU}]
If $\zt:E\to M$ is a vector bundle, then $\we^r\sT E$ is canonically a graded bundle of degree $r$ with respect to the projection
\begin{equation*}\we^r\sT\zt:\we^r\sT E\to \we^r\sT M\,.\end{equation*}
For $r=2$, the adapted coordinates on $\we^2E$ are $(x^\zr, y^a,{\dot x}^{\zm\zn}, y^{\zs b}, z^{cd} )$, ${\dot x}^{\zm\zn}=-{\dot x}^{\zn\zm}$, $z^{cd}=-z^{dc}$, coming from the decomposition of a bivector
\begin{equation*}\wedge ^2\sT E\ni u = \frac{1}{2} {\dot x}^{\zm\zn} \frac{\partial}{\partial x^\zm}\wedge \frac{\partial }{\partial x^\zn} + y^{\zs b}\frac{\partial }{\partial x^\zs}\wedge \frac{\partial }{\partial y^b} +\frac{1}{2} {z}^{cd} \frac{\partial }{\partial y^c}\wedge \frac{\partial }{\partial y^d}\,,
\end{equation*}
are of degrees $0,1,0,1,2$, respectively.
\end{Example}

Note that objects similar to graded bundles have been used in supergeometry by {\v{S}evera \cite{Severa:2005}, Voronov \cite{Voronov:2001qf}, Roytenberg \cite{Roytenberg:2001}} et al. under the name \emph{N-manifolds}. However, we will work with classical, purely even manifolds only.

Mimicking the definition of a vector bundle morphism, we get the following.
\begin{Definition} \textit{Morphisms} in the  \textit{category of graded bundles} are represented by commutative diagram of smooth maps
\begin{equation*}\xymatrix{
F^1\ar[rr]^{\Phi} \ar[d]^{\zt_1} && F^2\ar[d]^{{\zt_2}} \\
M_1\ar[rr]^{{\varphi}} && M_2 }
\end{equation*}
which are morphisms of graded spaces in fibers, i.e. which locally preserve the weight of homogeneous coordinates.
\end{Definition}

\begin{Example}
Any smooth map $\phi:M_1\to M_2$ induces a canonical morphism of graded bundles $\sT^k\phi:\sT^kM_1\to\sT^kM_2$.
\end{Example}

One can pick an atlas of $F$ consisting of charts for which the degrees of homogeneous local coordinates $(x^{A}, y_{w}^{a})$ are $\wu(x^{A}) =0$ and  $\wu(y_{w}^{a}) = w$, \ $1\leq w \leq k$, where $k$ is the degree of the graded bundle. The local changes of coordinates  are of the form
\begin{eqnarray}\label{eqn:translaws}
x'^{A} &=& x'^{A}(x),\\
\nonumber y'^{a}_{w} &=& y^{b}_{w} T_{b}^{\:\: a}(x) + \sum_{\stackrel{1<n  }{w_{1} + \cdots + w_{n} = w}} \frac{1}{n!}y'^{b_{1}}_{w_{1}} \cdots y'^{b_{n}}_{w_{n}}T_{b_{n} \cdots b_{1}}^{\:\:\: \:\:\:\:\:a}(x),
\end{eqnarray}
where $T_{b}^{\:\: a}$ are invertible and $T_{b_{n} \cdots b_{1}}^{\:\:\: \:\:\:\:\:a}$ are symmetric in indices $b_1,\dots,b_n$.

Note that the homogeneity structure in the typical fiber of a graded bundle $F$, i.e. the action $h:\R\ti\R^\bd\to\R^\bd$, is preserved under the transition functions, that defines a globally defined homogeneity structure $h:\R\ti F\to F$. In local homogeneous coordinates it reads
\begin{equation*}h_t(x^A,y_w^a)=(x^A,t^{w}y^a_w)\,.\end{equation*}
We call a function $f:F\to \R$  \textit{homogeneous of degree (weight) $w$} if
\begin{equation*}f\circ h_t=t^w\cdot f\,.\end{equation*}
The whole information about the degree of homogeneity is contained in the  \textit{weight vector field} (for vector bundles called the \textit{Euler vector field})
\begin{equation*}\nabla_F=\sum_swy^a_w\pa_{y^a_w}\,,\end{equation*}
so $f:F\to \R$ is homogeneous of degree $w$ if and only if $\nabla_F(f)=w\cdot f$.
Clearly, the fiber bundle morphism $\Phi$ is a smooth map which relates the weight vector fields $\nabla_{F^1}$ and $\nabla_{F^2}$.

Since, according to (\ref{eqn:translaws}), in the transition functions for a graded bundle, coordinates of weights not greater than $j$ depend only on coordinates of weights not greater than $j$, we can ``forget'' coordinates of higher weights, reducing a graded bundle $F$ to its part $F_j$ of degree $j$. The coordinate transformations for the canonical projection $F_j\to F_{j-1}$ are linear modulo a shift by a polynomial in variables not greater than $j$, so the fibrations $F_j\to F_{j-1}$ are affine. In this way we get for any graded bundle $F$ of degree $k$, like for jet bundles, a tower of affine fibrations
\begin{equation}\label{bu:tower}
F=F_{k} \stackrel{\tau^{k}}{\longrightarrow} F_{k-1} \stackrel{\tau^{k-1}}{\longrightarrow}   \cdots \stackrel{\tau^{3}}{\longrightarrow} F_{2} \stackrel{\tau^{2}}{\longrightarrow}F_{1} \stackrel{\tau^{1}}{\longrightarrow} F_{0} = M\,.
\end{equation}
In the case of the canonical graded bundle $F=\sT^kM$, we get exactly the tower of projections of jet bundles:
\begin{equation}\label{bu:towertk}
\sT^kM \stackrel{\tau^{k}}{\longrightarrow} T^{k-1}M \stackrel{\tau^{k-1}}{\longrightarrow}   \cdots \stackrel{\tau^{3}}{\longrightarrow} \sT^{2}M \stackrel{\tau^{2}}{\longrightarrow}\sT M \stackrel{\tau^{1}}{\longrightarrow} F_{0} = M\,.
\end{equation}

The fundamental fact is that graded bundles and homogeneity structures are actually equivalent concepts.

\begin{Theorem}[Grabowski-Rotkiewicz \cite{Grabowski:2012}]
For any homogeneity structure $h$ on a manifold $F$, there is a smooth submanifold $M$ of $F$, a non-negative integer $k\in\mathbb N$, and an $\R$-equivariant map
$\Phi_h^k:F\to \sT^kF_{|M}$
which identifies $F$ with a graded submanifold of the graded bundle $\sT^kF$. In particular, there is an atlas on $F$ consisting of local homogeneous coordinates.
\end{Theorem}

Since  {morphisms} of two homogeneity structures are defined as smooth maps $\Phi:F_1\to F_2$ intertwining the $\R$-actions: $\Phi\circ h^1_t=h^2_t\circ\Phi$, this describes also morphism of graded bundles. Consequently, a  \textit{graded subbundle} of a graded bundle $F$ is a smooth submanifold $S$ of $F$ which is invariant with respect to homotheties, $h_t(S)\subset S$ for all $t\in\R$.

The principle that says \textit{multiplication by reals is enough} has now the following consequences for vector bundles.

\begin{Corollary}
A smooth map $\Phi:E_1\to E_2$  between the total spaces of two vector bundles $\zp_i:E_i\to M_i$, $i=1,2$, is a morphism of vector bundles if and only if it intertwines the multiplications by reals:
\begin{equation*}\Phi(t\cdot v)=t\cdot\Phi(v)\,.\end{equation*}
In this case, the map $\phi=\Phi_{|M_1}$ is a smooth map between the base manifolds covered by $\Phi$.
\end{Corollary}

\begin{Corollary}
A submanifold $S$ of a vector bundle $\pi:E\to M$ is a vector subbundle if and only if it is invariant with respect to homotheties (multiplication by reals):
\begin{equation*}h_t(S)\subset S\,.\end{equation*}
In this case, $\zp(S)$ is a submanifold of M which is
the support of the vector subbundle $S$.
\end{Corollary}

\section{Double structures}\label{sec:DoubleStructures}

In geometry and applications one often encounters \textit{double vector bundles}, i.e.
manifolds equipped with two vector bundle structures which are  compatible in a categorical sense.   They were defined by Pradines  \cite{P} and studied by Mackenzie \cite{Mackenzie:2005}, Konieczna(Grabowska), and Urba\'nski \cite{KU} as  \textit{vector bundles in the category of vector bundles}. More precisely:

\begin{Definition}  A  \textit{double vector bundle} $(D;A,B;M)$ is a system of four
vector bundle structures
$$
\xymatrix@R+4mm @C+4mm{
  D \ar[d]_{\qDA} \ar[r]^{\qDB}
                & B \ar[d]^{\qB} &   \\
  A  \ar[r]^{\qA}
                & M         }
$$
in which $D$ has two vector bundles structures, on bases $A$ and
$B$.   The latter are themselves vector bundles on $M$, such that each of
the four structure maps of each  vector bundle structure on $D$
(namely the bundle projection, zero section, addition and scalar
multiplication) is a morphism of vector bundles with respect to the
other structures.
\end{Definition}

\medskip
In the above figure, we refer to $A$ and $B$ as the  \emph{side bundles} of
$D$, and to $M$ as the \emph{double base}. In the two side bundles, the
addition and scalar multiplication are denoted by the
usual symbols $+$ and juxtaposition, respectively.
We distinguish the two zero-sections, writing $0^A: M\lon A$, $m\mapsto 0^A_m$, and $0^B:
M\lon B$, $m\mapsto 0^B_m$.

In the vertical bundle structure on $D$ with base $A$, the vector bundle operations are denoted by $\Aplus$ and $\Atimes$, with $\zeroDA: A\lon D$, $a\mapsto \zeroDA_a$, for the zero-section. Similarly, in the horizontal bundle structure on $D$ with base $B$ we write $\Bplus$ and $\Btimes$, with $\zeroDB: B\lon D$, $b\mapsto \zeroDB_b$, for the zero-section.
%The two structures on $D$, namely $(D,\qDB,B)$ and $(D,\qDA,A)$ will also be denoted, respectively, by $\tilde{D}_B$ and $\tilde{D}_A$,  and called the \textit{horizontal bundle structure} and the \textit{vertical bundle structure}.

The condition that each vector bundle operation in $D$ is a morphism with respect to the other is equivalent to the following conditions, known as the \emph{interchange laws}:
\begin{eqnarray*}\label{Eqt:interchange1}
(d_1\Bplus d_2)\Aplus (d_3\Bplus d_4)&=& (d_1\Aplus d_3)\Bplus
(d_2\Aplus d_4),\\
\label{Eqt:interchange2} t\Atimes(d_1\Bplus
d_2)&=&t\Atimes d_1\Bplus t\Atimes d_2,\\\label{Eqt:interchange3}
t\Btimes(d_1\Aplus d_2)&=&t\Btimes d_1\Aplus t\Btimes
d_2,\\\label{Eqt:interchange4} t\Atimes(s\Btimes d)&=&
s\Btimes(t\Atimes d),\\\label{Eqt:interchange5}
\zeroDA_{a_1+a_2}&=&\zeroDA_{a_1}\Bplus \zeroDA_{a_2},\\
\zeroDA_{ta}&=&t\Btimes \zeroDA_{a},\\\label{Eqt:interchange6}
\zeroDB_{b_1+b_2}&=&\zeroDB_{b_1}\Aplus \zeroDA_{b_2},\\
\zeroDB_{tb}&=&t\Atimes \zeroDB_{b}.
\end{eqnarray*}

\medskip
We can extend the concept of a double vector bundle of Pradines to double graded bundles. However, thanks to our simple description in terms of a homogeneity structure, the `diagrammatic' definition of Pradines can be substantially simplified. As two graded bundle structure on the same manifold are just two homogeneity structures, the obvious concept of compatibility leads to the following:

\begin{Definition} A  \textit{double graded bundle} is a manifold equipped with two homogeneity structures $h^1,h^2$ which are  \textit{compatible} in the sense that
$$h^1_t\circ h^2_s=h^2_s\circ h^1_t\quad \text {for all\ } s,t\in\R\,.$$
\end{Definition}

The above condition can also be formulated as commutation of the corresponding weight vector fields, $[\nabla^1,\nabla^2]=0$. For vector bundles this is equivalent to the concept of a double vector bundle in the sense of Pradines.

\begin{Theorem}[Grabowski-Rotkiewicz \cite{GR0}]
The concept of a double vector bundle, understood as a particular double graded bundle in the above sense, coincides with that of {Pradines}.
\end{Theorem}

All this can be extended to  {$n$-fold graded bundles} in the obvious way:

\begin{Definition} A  {$n$-fold graded bundle} is a manifold equipped with $n$ homogeneity structures $h^1,\dots,h^n$ which are  {compatible} in the sense that
\vskip-.3cm$$h^i_t\circ h^j_s=h^j_s\circ h^i_t\quad \text {for all\ } s,t\in\R\quad \text{and}\quad i,j=1,\dots,n\,.$$
\end{Definition}

\begin{Example}
If $\zt:F\to M$ is a graded bundle of degree $k$, then there are  {canonical lifts} of the graded structure to the tangent and to the cotangent bundle.   In this way $\sT F$ and $\sT^*F$ carry canonical double graded bundle structure: one is the obvious vector bundle, the other is the lifted one (of degree $k$). There are also lifts of graded structures on $F$ to $\sT^r F$.

For local homogeneous coordinates $(x^A,y^a_w)$ on $F$, the adopted coordinates $(x^A,y^a_w,\dot x^B,\dot y^b_{w'})$ on $\sT F$ have degrees $0,w,0,w'$, respectively.
The coordinates in fibers of $\sT^*F$, dual to $\dot x^B,\dot y^b_{w'}$ have degrees $k$ and $k-w'$.

\medskip
In particular, if $\zt:E\to M$ is a vector bundle, then $\sT E$ and $\sT^*E$ are double vector bundles. There is a canonical isomorphism of double vector bundles
$$\sT^*E^*\simeq\sT^*E\,.$$
\end{Example}

\begin{Definition} A double graded bundle whose one structure is linear we will call a \emph{$\GrL$-bundle}.
\end{Definition}
Canonical examples of $\GrL$-bundles are $\sT F$ and $\sT^*F$ with their lifted and linear structures. Another canonical examples are tensor bundles.

\begin{Example}
if $\zt:E\to M$ is a vector bundle, then $\we^k\sT E$ is canonically a $\GrL$-bundle:
$${\xymatrix@R-2mm @C-2mm{ & \wedge^k\sT E \ar[ld] \ar[rd] & \cr
        \quad E \ar[rd] & & \wedge^k\sT M \ar[ld]  \cr & M  & }}\,.
$$
\end{Example}

\section{Linearisation of graded bundles}\label{sec:Linearisation}

The possibility of constructing mechanics on graded bundles is based on the following generalization of the embedding $\sT^kQ\hookrightarrow\sT\sT^{k-1}Q$ introduced in \cite{Bruce:2014}.  However, we will not discuss the applications of graded bundles in geometric mechanics and stick to pure mathematical aspects.

\begin{Theorem}[Bruce-Grabowska-Grabowski \cite{Bruce:2016}]
There is a canonical  \textit{linearisation functor} $\operatorname{l}:\GrB\to\GrL$ from the category of graded bundles into the category of $\GrL$-bundles   which assigns, for an arbitrary graded bundle $F_k$ of degree $k$, a canonical $\GrL$-bundle $\operatorname{l}(F_k)$ of bi-degree $(k-1,1)$ which is linear over $F_{k-1}$,   called the {linearisation of $F_k$},   together with a  {graded} embedding $\zi:F_k\hookrightarrow \operatorname{l}(F_k)$ of $F_k$ as an affine subbundle of the vector bundle $\operatorname{l}(F_k)\to F_{k-1}$.
\end{Theorem}

Elements of $F_k\subset \operatorname{l}(F_k)$ may be viewed as  {`holonomic vectors'} in the linear-graded bundle $\operatorname{l}(F_k)$. We have $\pLinr(\sT^kM)\simeq\sT\sT^{k-1}M$ and
\begin{equation}\label{lin:tk}\zi:\sT^kM\hookrightarrow \operatorname{l}(\sT^kM)\simeq\sT\sT^{k-1}M\end{equation}
is the canonical embedding of $\sT^kM$ as holonomic vectors in $\sT\sT^{k-1}M$.

For instance,
if $(x^a, y^A,z^j)$ are coordinates on a graded bundle $F_2$ of degrees $0,1,2$, respectively, then, the induced coordinate system on $\operatorname{l}(F_2)$ is
$$(x^a,y^A,\dot y^{B},\dot z^j)\,,$$
where $x^a$, $y^A$, $\dot y^B$, and $\dot z^j$ are of bi-degree $(0,0)$, $(1,0)$, $(0,1)$, and $(1,1)$, respectively.

The transformation laws for the extra coordinates  are obtained by differentiation;
\begin{align*}
&\dot{y}^{A} = \dot{y}^{B}T_{B}^{A}(x),&\\
&\dot{z}^{i} = \dot{z}^{j}T_{j}^{\:\:i}(x)+ \dot{y}^{B}y^{A}T_{AB}^{i}(x)\,.&
\end{align*}
Thus,
$$(x^a,y^A,\dot y^{B},\dot z^j)\mapsto (x^a,y^A)$$
is a linear fibration over $F_{1}$.
The embedding $\zi:F_2\hookrightarrow \operatorname{l}(F_2)$ reads
$$\zi(x^a,y^A,z^j)=(x^a,y^A,y^{A},2z^j)\,.
$$

\begin{Example}\label{exm:6} For a Lie groupoid $G \rightrightarrows M$, consider the subbundle $\sT^{k}G^{\underline{s}} \subset \sT^{k}G$ consisting of all higher order velocities tangent to source-leaves. The bundle
\begin{equation*}
 F_{k} = {A}^{k}(G) := \left. \sT^{k}G^{\underline{s}}\, \right|_{M},
\end{equation*}
inherits graded  bundle structure of degree $k$ as a graded subbundle of $\sT^{k}G$. Of course, $A=A^1(G)$ can be identified with the Lie algebroid of $G$.

\begin{Theorem}[Bruce-Grabowska-Grabowski \cite{Bruce:2016}]\label{thm:linearisation}
The linearisation of ${A}^{k}(G)$ is given as
\begin{equation*}
\pLinr({A}^{k}(G) ) \simeq \{ (Y,Z) \in {A}(G) \times \sT{A}^{k-1}(G)|\quad {\rho}(Y) = \sT {\tau}(Z)  \}\,,
\end{equation*}
viewed as a vector bundle over ${A}^{k-1}(G)$ with respect to the obvious projection of part $Z$ onto ${{A}^{k-1}(G)}$,   where ${\rho} : {A}(G) \rightarrow \sT M$ is the standard anchor of the Lie algebroid and ${\tau}: {A}^{k-1} (G) \rightarrow M$ is the obvious projection.
\end{Theorem}
Let us remark that the above linearization is canonically a weighted Lie algebroid, a  {Lie algebroid prolongation} in the sense of  Cari\~{n}ena \&  Mart\'{\i}nez \cite{Carinena:2001} and Popescu \& Popescu \cite{Popescu:2001}. The corresponding groupoid prolongations have been described by Saunders \cite{Saunders:2004}.
\end{Example}

Applying the linearisation functor consecutively to a graded bundle of degree $k$, we arrive at a $k$-fold graded bundle od degree $(1,\dots,1)$,   i.e. at a $k$-fold vector bundle. This functor from $\GrB[k]$ to $\VB[k]$ we call a \emph{total linearisation}. Its range consists of $k$-fold vector bundles equipped with an action of the symmetry group $S_k$ permuting the order of vector bundle structures (\emph{symmetric $k$-fold vector bundles}).

\begin{Theorem}[Bruce-Grabowski-Rotkiewicz \cite{Bruce:2016a}]
There is a canonical functor $\Linr[k]:\GrB[k]\to\VB[k]$ from the category of graded bundles of degree $k$ into the category of $k$-fold vector bundles.   It gives an equivalence of $\GrB[k]$ with the subcategory (not full) $\SymmVB[k]$ of {symmetric $k$-fold vector bundles}. There is a canonical {graded} embedding $\zi[k]:F_k\hookrightarrow \operatorname{L}(F_k)$ of $F_k$ as a subbundle of  {symmetric (holonomic) vectors}.
\end{Theorem}

\begin{Example} We have $\Linr(\sT^kM)\simeq\sT^{(k)}M$, where $\sT^{(k)}M=\sT\sT\cdots\sT M$ is the iterated tangent bundle. The action of $S_k$ comes from iterations of the canonical ``flips" $\kappa:\sT\sT M\to\sT\sT M$.
\end{Example}

\begin{Example}
Continuing Example \ref{exm:6} and using Theorem \ref{thm:linearisation}, we get that $\Linr(A^{3}(G))$ equals
$$
\left\{(X,Y,Z) \in A(G)\times \sT A(G) \times \sT^{(2)}A(G)~ |~\rho(X) = \sT \pi (Y)\,,\  \sT \rho(Y) = \sT^{(2)}\pi(Z)  \right\}\,,
$$
where $\sT^{(l)} =  \sT \sT \cdots \sT$ ($l$-times),  $\pi : A(G) \rightarrow M$ is the standard projection, and $\rho :A(G)\rightarrow \sT M$  is the anchor of the Lie algebroid.  To see this we note that the linearisation functor as a subfunctor of the tangent functor respects products and commutes with the tangent functor. In particular, we have
$$\pLinr^{(2)}(A^{3}(G)) \subset \pLinr(A(G)\times \sT A^{2}(G)) = A(G)\times \sT \pLinr(A^{2}(G)),$$
Thus
$$\Linr(A^{3}(G)) \subset A(G) \times \sT A(G) \times \sT^{(2)}A(G).$$
Then, applying the tangent functor to the condition given in Theorem \ref{thm:linearisation}, we arrive at the desired result.
\end{Example}

The previous example generalises directly to the higher order case just by iteration. We this have the following theorem.
\begin{Theorem}
The full linearisation of ${A}^{k}(G)$ is given as
\begin{align*}&\Linr(A^{k}(G)) =
\left\{(X_{1},X_{2}, \cdots ,X_{k}) \in A(G)\times \sT A(G) \times \sT^{(2)}A(G) \cdots \times \sT^{(k-1)}A(G)  | \right. \\
& \left. \rho(X_{1}) = \sT \pi (X_{2}), \sT \rho(X_{2}) = \sT^{(2)} \pi (X_{3}), ~\cdots, ~ \sT^{(k-2)} \rho(X_{k-1}) = \sT^{(k-1)}\pi(X_{k})  \right\},
\end{align*}
where $\sT^{(l)} =  \sT \sT \cdots \sT$ ($l$-times),  $\pi : A(G) \rightarrow M$ is the standard projection, and $\rho :A(G)\rightarrow \sT M$  is the anchor of the Lie algebroid.
\end{Theorem}

\section{Weighted Lie groupoids and algebroids}\label{sec:WeightedLie}

Besides the compatibility of graded bundle structures, we can consider a compatibility of a graded bundle structure with some other geometric structures, e.g. a Lie algebroid or a Lie groupoid structure. Thanks to the fact that a graded bundle structure can be expressed in terms of an $(\R,\cdot)$-action, there is an obvious natural concept of such a compatibility.

\begin{Definition} A  \textit{weighted Lie groupoid} (resp., a  \textit{weighted Lie algebroid}) of degree $k$ is a Lie groupoid (resp., Lie algebroid) equipped with a homogeneity structure $h$ of degree $k$ such that homotheties $h_t$ act as Lie groupoid (resp., Lie algebroid) morphisms.
\end{Definition}

We use the name `weighted', as the term  \textit{graded Lie algebroids} is already used in various meanings. Note that weighted Lie groupoids (algebroids) of degree $1$ have already appeared in the literature under the name \emph{$\VB$-groupoids} (\emph{$\VB$-algebroids}).

\begin{Example} If $\mathcal{G}$ is a Lie groupoid (algebroid), then $\sT^k\mathcal{G}$ is canonically a weighted Lie groupoid (algebroid) of degree $k$.
\end{Example}

\medskip
Note that the compatibility condition between the extra homogeneity structure on $\mathcal{G}$ and its Lie algebroid structure we use in applications for mechanics   is that the double vector bundle morphism associated with the Lie algebroid structure
$\epsilon : \sT^{\ast} \mathcal{G} \simeq \sT^{\ast}\mathcal{G}^{\ast}\rightarrow \sT \mathcal{G}^{\ast}$
is a  morphisms of triple graded bundles.

\begin{Theorem}[Bruce-Grabowska-Grabowski \cite{Bruce:2014a}] There is a one-to-one correspondence between weighted Lie groupoids of degree $k$ with simple-connected source fibers and  {integrable} weighted Lie algebroids of degree $k$. In other words, homogeneity structures compatible with Lie structures can be differentiated and integrated.
\end{Theorem}

\begin{Example} Let $G$ be a Lie groupoid with the Lie algebroid {{$\mathcal{G}$}}.  The weighted Lie algebroid for $\sT^k{G}$ is $\sT^k\mathcal{G}$.
\end{Example}
%\section*{Acknowledgments}
%This work was supported by the  Polish National Science Centre grant under the contract number DEC-2012/06/A/ST1/00256.

\vskip1cm

\noindent Andrew James Bruce\\
\emph{Institute of Mathematics, Polish Academy of Sciences,}\\ {\small \'Sniadeckich 8,  00-656 Warszawa, Poland}\\ {\tt andrew@impan.pl}\\

\noindent Katarzyna Grabowska\\
\emph{Faculty of Physics,
                University of Warsaw} \\
               {\small Pasteura 5, 02-093 Warszawa, Poland} \\
                 {\tt konieczn@fuw.edu.pl} \\

 \noindent Janusz Grabowski\\\emph{Institute of Mathematics, Polish Academy of Sciences}\\{\small \'Sniadeckich 8, 00-656 Warszawa,
Poland}\\{\tt jagrab@impan.pl}

\end{document}